\documentclass[a4paper,10pt,reqno]{amsart}
\usepackage[a4paper,tmargin=40mm,bmargin=35mm,hmargin=30mm]{geometry}
\usepackage[T1]{fontenc}
\usepackage{lmodern}
\usepackage{amsmath}
\usepackage{amssymb}
\usepackage{amsthm}
\usepackage{stmaryrd}
\usepackage{tikz}
\usepackage[all]{xy}
\usepackage{booktabs}
\usepackage{multirow}
\usepackage{multicol}
\usepackage{here}
\usepackage{aliascnt}
\usepackage{hyperref}
\usepackage[noabbrev]{cleveref}

\newcommand{\NewTheorem}[2]{
	\newaliascnt{#1}{TheoremEnvironment}
	\newtheorem{#1}[#1]{#1}
	\aliascntresetthe{#1}
	\crefname{#1}{#1}{#2}
	\Crefname{#1}{#1}{#2}
}

\theoremstyle{definition}

\NewTheorem{Definition}{Definitions}
\NewTheorem{Remark}{Remarks}
\NewTheorem{Axiom}{Axioms}
\NewTheorem{Example}{Examples}
\NewTheorem{Observation}{Observations}
\NewTheorem{Convention}{Conventions}
\NewTheorem{Notation}{Notations}
\NewTheorem{Setting}{Settings}
\NewTheorem{Question}{Questions}
\NewTheorem{Answer}{Answers}
\NewTheorem{Conjecture}{Conjectures}
\NewTheorem{Problem}{Problems}
\NewTheorem{Solution}{Solutions}
\NewTheorem{Goal}{Goals}
\NewTheorem{Comment}{Comments}
\NewTheorem{Aim}{Aims}
\NewTheorem{Caution}{Cautions}
\NewTheorem{Exercise}{Exercises}
\NewTheorem{Examples}{Examples}
\theoremstyle{plain}
\NewTheorem{Proposition}{Propositions}
\NewTheorem{Lemma}{Lemmas}
\NewTheorem{Theorem}{Theorems}
\NewTheorem{Corollary}{Corollaries}

\crefname{enumi}{}{}
\Crefname{enumi}{}{}
\creflabelformat{enumi}{(#2#1#3)}
\crefname{enumii}{}{}
\Crefname{enumii}{}{}
\creflabelformat{enumii}{(#2#1#3)}
\crefname{enumiii}{}{}
\Crefname{enumiii}{}{}
\creflabelformat{enumiii}{(#2#1#3)}

\makeatletter
\renewcommand{\p@enumii}{}
\renewcommand{\p@enumiii}{}
\makeatother

\numberwithin{equation}{section}
\crefname{equation}{}{}
\Crefname{equation}{}{}
\creflabelformat{equation}{(#2#1#3)}

\newcommand{\SwapSymbols}[1]{
	\expandafter\let\expandafter\temporarysymbol\csname #1\endcsname
	\expandafter\let\csname #1\expandafter\endcsname\csname var#1\endcsname
	\expandafter\let\csname var#1\endcsname\temporarysymbol
}

\SwapSymbols{epsilon}
\SwapSymbols{phi}
\SwapSymbols{Gamma}
\SwapSymbols{Delta}
\SwapSymbols{Theta}
\SwapSymbols{Lambda}
\SwapSymbols{Xi}
\SwapSymbols{Pi}
\SwapSymbols{Sigma}
\SwapSymbols{Upsilon}
\SwapSymbols{Phi}
\SwapSymbols{Psi}
\SwapSymbols{Omega}

\newcommand{\bbM}{\mathbb{M}}
\newcommand{\bbN}{\mathbb{N}}

\newcommand{\cF}{\mathcal{F}}
\newcommand{\cG}{\mathcal{G}}

\newcommand{\cI}{\mathcal{I}}

\newcommand{\cX}{\mathcal{X}}
\newcommand{\cY}{\mathcal{Y}}

\newcommand{\To}{\longrightarrow}

\DeclareMathOperator{\Hom}{Hom}
\DeclareMathOperator{\End}{End}
\DeclareMathOperator{\Ext}{Ext}
\DeclareMathOperator{\Tor}{Tor}
\DeclareMathOperator{\id}{id}

\DeclareMathOperator{\Gid}{Gid}
\DeclareMathOperator{\RHom}{{\bf R}Hom}

\DeclareMathOperator{\Ann}{Ann}

\DeclareMathOperator{\Ker}{Ker}
\DeclareMathOperator{\Coker}{Coker}
\let\Im\relax
\DeclareMathOperator{\Im}{Im}

\DeclareMathOperator{\Spec}{Spec}

\DeclareMathOperator{\Ass}{Ass}
\DeclareMathOperator{\Gpd}{Gpd}

\DeclareMathOperator{\height}{ht}
\DeclareMathOperator{\depth}{depth}

\title{Bass numbers and endomorphism rings of Gorenstein injective modules}

\author{Reza Sazeedeh}

\address{Department of Mathematics, Urmia University, P.O.Box: 165, Urmia, Iran}
\email{rsazeedeh@ipm.ir and r.sazeedeh@urmia.ac.ir}

\begin{document}

\begin{abstract}
Let $R$ be a commutative noetherian ring admitting a dualizing complex and let $\mathfrak p$ be a prime ideal of $R$. In this paper we investigate when $G(R/\frak p)$ is an $R_{\frak p}$-module. We give some necessary and sufficient conditions under which $G(R/\frak p)$ is an $R_{\frak p}$-module. We also study the Bass numbers of $G(R/\frak p)$ and we show that if $\Gid_RR/\frak p$ is finite, then $\mu^i(\frak q,G(R/\frak p))$ is finite for all $i\geq 0$ and all $\frak q\in\Spec R$. If $\Gpd_RR/\frak p$ is finite, then $\mu^i(\frak p,G(R/\frak p))$ is finite for all $i\geq 0$.  We define a subring $S(\frak p)_{\frak p}$ of $\End_{R_{\frak p}}(G(R_{\frak p}/\frak pR_{\frak p}))$ and we show that it is noetherian and contains a subring which is a quotient of $\widehat{R_{\frak p}}$.  
\end{abstract}

\maketitle
\tableofcontents

\section{Introduction}
Throughout this paper, assume that $R$ is a commutative noetherian ring with an identity $1_R$ admitting a dualizing complex $D$ and assume that $\frak p$ is a prime ideal of $R$. With this conditions, Enochs and Iacob [EI] showed that the class of Gorenstein injective modules is enveloping. Our main aim of this paper is to study the Bass numbers and endemorphism ring of  $G(R/\frak p)$, where $G(R/\frak p)$ denotes the Gorenstein injcetive envelope of $R/\frak p$.  

 The motivation of this work originally goes back to a work of Matlis [M] in 1958. He showed that $\End_R(E(R/\frak p))$ is isomorphic in a natural way to $\widehat{R_{\frak p}}$, the completion of $R_{\frak p}$ with respect to the maximal ideal $\frak pR_{\frak p}$. In this paper we want to find out whether this property can be extended to $\End_R(G(R/\frak p))$, where $G(R/\frak p)$ is the injective envelope of $R/\frak p$. The first attempt in this material was done by Belshoff and Enochs [BE] for $G(R/\frak m)$ where $(R,\frak m)$ was a Gorenstein local ring. To be more precise, they introduced a subring  $S$ of $\End_R(G(R/\frak m))$ and they proved that $S\cong\widehat{R}$ when $R$ is a Gorenstein local ring of dimension one. In the end of the paper, they posed some questions about the ring $S$ which is quietly unknown. For example when $S$ is commutative or is noetherian  and what can say about $G(R/\frak m)$ as a left $S$-module. In this paper, for any prime ideal $\frak p$ of a commutative noetherian ring $R$, we define the subring $S(\frak p)$ of $\End_R(G(R/\frak p))$ and  we try to answer to the above questions.

 In Section 2, we study the Gorenstein injective envelopes of modules. Some of the results are well known for modules over Gorenstein rings and we prove them without the condition Gorenstein on $R$. for the new results, assume that $\frak a$ is an ideal of $R$. It is a well known result that if $M$ is an $R$-module such that $\Gamma_{\frak a}(M)=0$, then $\Gamma_{\frak a}(E(M))=0$ where $E(M)$ is the injective envelope of $M$. In \cref{exx}, we show that $\Gamma_{\frak a}(G(M))$ is not necessarily zero, where $G(M)$ is the Gorenstein injective envelope of $M$. Further, in \cref{cap}, we show that $\Gamma_{\frak a}(G(M))$ is reduced (i.e. it does not contain an injective submodule). In \cref{evg}, we show that every injective extension of a module $M$ can be extended to a Gorenstein injective preenvelope of $M$. We show that if $M$ is $\frak a$-torsion (i.e. $\Gamma_{\frak a}(M)=M)$, then so is $G(M)$. At the end of this section, we show that every Gorenstein injective module is the direct sum of  an injective module and a reduced Gorenstein injective module (see \cref{go}).

  In Section 3, we first show that $G(R/\frak p)$ is indecomposable. It is know that for  every indecmoposable injective module $E$, there exists a prime ideal $\frak p$ of $R$ such that $E=E(R/\frak p)$. In \cref{ex,eexx}, we show that this result can not be generalized for indecmoposable Gorenstein injective modules. Assume that $\sup D=s$ and  for any integer $k$ we set $X_k(\frak p)=\{\frak q\in V(\frak p)|\hspace{0.1cm} {\rm ht}\frak p-\sup D_{\frak p}=k\}$. One of the main result is the following theorem.
  
  \begin{Theorem}
Let $\frak p\in\Spec R$ such that $\height \frak p-\sup D_{\frak p}=t$. Then there exists a finite filtration $0\subset G_d\subset G_{d-1}\subset\dots \subset G_t=G(R/\frak p)$ of Gorenstein injective submodules of $G(R/\frak p)$ such that $G_k/G_{k+1}\cong\bigoplus_{\frak q\in X_{k-s}(\frak p)}\Tor_{k-s}^R(E(R/\frak q),\RHom_R(D,G(R/\frak p))$ is reduced for each $t\leq k\leq d$. In particular, 
$$G_t/G_{t+1}\cong \Tor_t^R(E(R/\frak p),\RHom_R(D,G(R/\frak p)))\cong G(R/\frak p)_{\frak p}.$$
\end{Theorem}
As a conclusion of the above theorem  we have $G(R/\frak p)_{\frak p}\cong G(k(\frak p))$ where $k(\frak p)\cong R_{\frak p}/\frak pR_{\frak p}$. It is well known that $E(R/\frak p)$ is in a natural way an $R_{\frak p}$-module and hence a natural question is given rise to ask whether $G(R/\frak p)$ is an  $R_{\frak p}$-module as well. The following theorem ensures when this occurs for a prime ideal $\frak p$ of $R$. 

\begin{Theorem}
Let $\frak p\in\Spec R$ with ${\rm ht}\frak p-\sup D_{\frak p}=t$. Then the following conditions are equivalent.

${\rm (i)}$ $\Ass_R G(R/\frak p)=\{\frak p\}$.

${\rm (ii)}$ $G(R/\frak p)$ is an $R_{\frak p}$-module.

${\rm (iii)}$ $G(R/\frak p)=G(k(\frak p))$ where $k(\frak p)=R_{\frak p}/\frak pR_{\frak p}$.

${\rm (iv)}$ $\Tor_i^R(E,\RHom_R(D,G(R/\frak p)))=0$ for all $i>t$.
\end{Theorem}
 If $\frak m$ is a maximal ideal, then $G(R/\frak m)$ is $\frak m$-torsion and so the above theorem implies that $G(R/\frak m)$ is an $R_{\frak m}$ module. We give some examples which show that $G(R/\frak p)$ is not an  $R_{\frak p}$-module in general if $\frak p$ is not maximal.
  
  In Section 4,  we study the Bass numbers of $G(R/\frak p)$ which have a key role in study of $\End(G(R/\frak p)$. We prove the following theorems. 

\begin{Theorem}
Let $\frak p\in\Spec R$ such that $\Gid_RR/\frak p$ is finite. Then $\mu^i(\frak q,G(R/\frak p))$ is finite for all $i\geq 0$ and all $\frak q\in\Spec R$.
\end{Theorem}
As a conclusion of this theorem, over a Gorenstein ring of finite Krull dimension $\mu^i(\frak q,G(R/\frak p))$ is finite for all $i\geq 0$ and all prime ideals $\frak q$ of $R$. Moreover, we have the following theorem.

\begin{Theorem}
Let $\frak p\in\Spec R$ such that $\Gpd_RR/\frak p$ is finite. Then $\mu^i(\frak p,G(R/\frak p))$ is finite for each $i\geq 0$. 
\end{Theorem}

In Section 5, we study the indemorphism rings of indecomposable Gorenstein modules in correspondence with $\End_R(R/\frak p)$. We  define a subring $S(\frak p)=\{f\in{\rm End}_R(G(R/\frak p))|\hspace{0.1cm} f(R/\frak p)\subset R/\frak p\}$ os $\End_R(G(R/\frak p))$ and we prove that $S(\frak p)$ is a local ring containing $R/\frak p$ with the maximal ideal $\frak n=\{f\in S(\frak p)|\hspace{0.1cm} f(R/\frak p)=0\}$. We set $S(\frak p)_{\frak p}=\{f_{\frak p}\in\End_{R_{\frak p}}(G(k(\frak p))|\hspace{0.1cm}f\in S(\frak p)\}$ and we show that $S(\frak p)_{\frak p}$ is local as well. We prove the following theorem. 
\begin{Theorem}
Let $\frak p\in\Spec R$ such that $\mu(\frak p,G(R/\frak p))$ is finite. Then $S(\frak p)_{\frak p}$ is noetherian.
\end{Theorem}
As an immediate consequence, if $R$ is a Gorenstein ring of Krull dimension $d$, then $S(\frak p)_{\frak p}$ is noetherian all $\frak p\in\Spec R$. We also show that  $\widehat{R_{\frak p}}/\Ann_{\widehat{R_{\frak p}}}G(R/\frak p)_{\frak p}$  is a (commutativ) subring of $S(\frak p)_{\frak p}$. In particular, if $R_{\frak p}$ is a Gorenstein local ring of dimension one, then $S(\frak p)_{\frak p}\cong\widehat{R_{\frak p}}$
For the basic notions about local cohomology modules, derived category and Gorenstein injective modules, we refer the reader to the text books [BS, EJ2,H].
\section{Gorenstein injective envelope of modules}
We start this section with some definitions which are used in throughout this paper.
\begin{Definition}
Let $\cF$ be class of $R$-module and $M$ be an $R$-module. An $\cF$-{\it preenvelope}  of $M$ is meant a morphism $\phi:M\To F$ such that $F\in\cF$ and for any morphism $f:M\To F'$ with $F'\in\cF$, there is a morphism $g:F\To F'$ such that $g\phi=f$. Furthermore, $\phi$ is an $\cF$-{\it envelope} if every $g:F\To F$ with $g\phi=g$ forces that $g$ is authomorphism.
\end{Definition}

\medskip

\begin{Definition} By [EJ2], an $R$-module $G$ is Gorenstein injective if there exists an exact sequence of injective $R$-modules 
$\cI:\dots\To E_1\To E_0\To E_{-1}\To \dots$ such that  $\Hom_R(E,\cI)$ is exact for all injective $R$-modules $E$ and $G=\Ker(E_0\To E_{-1})$. The exact sequence above is called a {\it complete} injective resolution of $G$. The class of Gorenstein injective modules is denoted by $\cG\cI$. Moreover $^{\bot}\cG\cI=\{X\in R\textnormal{-Mod}|\hspace{0.1cm} \Ext^1_R(X,G)=0 \hspace{0.1cm}\textnormal{for all Gorenstein injective modules}\hspace{0.1cm} G\}$. 
\end{Definition}

\medskip

 \begin{Definition} For any $R$-module  $M$, the {\it Gorenstein injective dimension} of $M$, denoted by $\Gid_RM$ is the least non-negative integer $n$ such that there is an exact sequence of $R$-modules  
 $$0\To M\To G^0\To G^1\To \dots \To G^n\To 0$$  
   with each $G_i$ Gorenstein injective. If there is no such $n$, we say that $\Gid_RM=\infty$. 
  The {\it Gorenstein projective dimension} of $M$, denoted by ${\rm \Gpd}_RM$ is defined dually. 
  \end{Definition}

  \medskip

\begin{Definition} A complex $D\in D_b^f(R)$ is {\it dualizing} for $R$ if it has finite injective dimension and the natural homothety morphism $\cX_D^R:R\To \RHom_R(D,D)$ is an isomorphism in $D(R)$. If $(R,\frak m)$ is a local ring, then the dualizing complex $D$ is {\it normalized} if $\sup D=\dim R.$ 
 \end{Definition}

  \medskip
\begin{Lemma}\label{sub}
Let $R$ be a noetherian ring and let $M\in^{\bot}{\cG\cI}$ be a submodule of some Gorenstein injective module $G$. Then $G(M)=E(M)$ and it is a direct summand of $G$. In particular $\widetilde{\cG\cI}\cap ^{\bot}\cG\cI=\widetilde{\cI}$, where $\widetilde{\cG\cI}$ ($\widetilde{\cI}$) is the class of all modules of finite (Gorenstein) injective dimension.   
\end{Lemma}
\begin{proof}
There exists an exact sequence of $R$-modules $0\To M\stackrel{\epsilon}\To G(M)\To X\To 0$ in which $G(M)$ is the Gorenstein injective envelope of $M$. On the other hand, by [K, Theorem 7.12], there exists a Gorenstein injective preenvelope of $M$ fitting into an exact sequence of modules $0\To M\To H \To Y\To 0$ such that $H$ is Gorenstein injective and $Y\in^{\bot}\cG\cI$. Since $G(M)$ is a direct summand of $H$, clearly $X$ is a direct summand of $Y$; and hence $X\in^{\bot}\cG\cI$.
Since $^{\bot}\cG\cI$ is closed under extensions, $G(M)\in\cG\cI\cap^{\bot}\cG\cI=\cI$, where $\cI$ denots the class of injective modules. Hence $G(M)$ is injective and clearly, it is the injective envelope of $M$. Since $G(M)$ is the Gorenstein injective envelope of $M$, there exists $\alpha:G(M)\To G$ such that $\alpha\epsilon=i$, where $i:M\To G$ is the inclusion morphism. On the other hand, since $G(M)$ is injective, there exists $\beta:G\To G(M)$ such that $\beta i=\epsilon$ and hence $\beta\alpha\epsilon=\epsilon$ so that $\beta\alpha$ is automorphism. This forces that $G(M)$ is a direct summand of $G$. For the second claim, if $M\in\widetilde{\cG\cI}\cap ^{\bot}\cG\cI$ with Gid$M=t$, by [EI, Corollary 2], the module $M$ has the Gorenstein injective envelope as $R$ admits a dualizing complex. Thus the first assertion implies that $G(M)=E(M)$ and an easy induction on $t$ gives the result.  
\end{proof}

\medskip

An $R$-module $M$ is said to be {\it reduced} if it has no a non-zero injective submodule. We have an immediate consequence.

  \begin{Corollary}\label{ccoo}
  Let $G$ be a Gorenstein injective $R$-module. Then $G$ is reduced if and only if it has no nontrivial submodules in $^{\bot}{\cG\cI}$. 
  \end{Corollary}
 \begin{proof}
 The proof is clear by \cref{sub} and noting the fact that $\cI\subseteq^{\bot}{\cG\cI}$.
 \end{proof} 
\medskip
 Enochs and Iacob [EI, Corollary 2] show that over a commutative noetherian ring admitting a dualizing complex, the class of Gorenstein injective modules is enveloping. Then we have the following proposition.
  
\begin{Proposition}\label{prred}
Let $M$ be an $R$-module and $M\To G$ be a Gorenstein injective preenvelope of $M$. Then there exists an injective $R$-module $E$ such that $G\cong G(M)\oplus E$.
\end{Proposition}
\begin{proof}
There exists the following commutative diagram with exact rows and $X\in ^{\bot}\cG\cI$
$$\xymatrix{0\ar[r]&M\ar@{=}[d]\ar[r]& G(M)\ar[d]^{f_1}\ar[r]&Y\ar[d]^{g_1} \ar[r]& 0\\
0\ar[r]&M\ar@{=}[d]\ar[r]&G \ar[r]\ar[d]^{f_2}& X\ar[r]\ar[d]^{g_2}& 0\\
0\ar[r]&M\ar[r]& G(M)\ar[r]&Y \ar[r]& 0}$$
and since $G(M)$ is Gorenstein envelope of $M$, the maps $f_2f_1$ and $g_2g_1$ are automorphisms and so $G=G(M)\oplus G/G(M)$ and $X=Y\oplus X/Y$. On the other hand, the diagram implies that $G/G(M)\cong X/Y\in ^{\bot}\cG\cI$ and so $E=G/G(M)$ is injective by \cref{sub}.  
\end{proof}

Let $N$ be a submodule of $M$. We recall from [EJ2] that $M$ is a Gorenstein essential extension  of $N$ if for every non-zero submodule $K$ of $M$ such that $K\in^{\bot}\cG\cI$, we have $K\cap N\neq 0$. 
\medskip
\begin{Proposition}\label{ess}
Let $M$ be an $R$-module. Then $G(M)$ is a Gorenstein essential extension of $M$.  
\end{Proposition}
\begin{proof}
Let $N$ be a non-zero submodule of $G(M)$ in $^{\bot}\cG\cI$.  It follows from \cref{sub} that $E(N)$ is a submodule of $G(M)$. Now, assume that $N\cap M=0$. Then $E(N)\cap M=0$ and so assume that $H$ is athe maximal Gorenstein injective submodule of $G(M)$ such that $E(N)\cap H=0$. If $E(N)+H$ is a proper submodule of $G(M)$, then $E(N)\cong (E(N)+H)/H\leq G(M)/H$ and so there exists a submodule $H_1$ of $G(M)$ containing $H$ such that $E(N)\oplus H_1/H=G(M)/H$. This implies that $E(N)+H_1=G(M)$ and since $E(N)\cap H=0$, the modular law implies that $E(N)\cap H_1=0$ which is a contradiction as $H$ is maximal. Therefore, $G(M)=E(N)\oplus H$; and hence $G(M)/M\cong H/M\oplus E(N)$ so that $H/M\in^{\bot}\cG\cI$. Thus $H$ is a Gorenstein injective preenvelope of $M$ so that $G(M)$ is a direct summand of $H$. Thus $E(N)=0$ which is a contradiction.  
\end{proof}

\medskip

\begin{Corollary}
Let $M$ be a an $R$-module and $E$ be an injective submodule of $G(M)$. Then $\Ass_RE\subset \Ass_RM$. 
\end{Corollary}
\begin{proof}
For any $\frak p\in\Ass_RE$, it follows from \cref{ess} that $E(R/\frak p)\cap M\neq 0$ and so $\frak p\in\Ass_RM$.
\end{proof}

\begin{Lemma}\label{secgor}
Let $\frak a$ be an ideal of $R$ and $G$ be a Gorenstein injective $R$-module. Then $H_{\frak a}^i(G)=0$ for all $i>0$. Furthermore if $R$ has a dualizing complex, then $\Gamma_{\frak a}(G)$ is Gorenstein injective.
\end{Lemma}
\begin{proof}
 See  [S2, Lemma 2.3].
\end{proof}

\medskip
\begin{Proposition}\label{cap}
Let $\frak a$ be an ideal of $R$ and let $M$ be an $R$-module such that $\Gamma_{\frak a}(M)=0$. Then $\Gamma_{\frak a}(G(M))$ is reduced.
\end{Proposition}
\begin{proof}
Given a submodule $X$ of $\Gamma_{\frak a}(G(M))$ in $^{\bot}\cG\cI$, it follows from \cref{ess} that $Y=X\cap M\neq 0$. We observe that $Y=\Gamma_{\frak a}(Y)\leq \Gamma_{\frak a}(M)=0$ so that $Y=0$ which is a contradiction. Now, \cref{secgor} and \cref{ccoo} imply that $\Gamma_{\frak a}(G(M))$ is reduced.
\end{proof}

The following example shows that the condition $\Gamma_{\frak a}(M)=0$ does not force to be $\Gamma_{\frak a}(G(M))=0$.

\medskip
\begin{Example}\label{exx}
Let $k$ be a field and set $R=k\llbracket X,Y\rrbracket/(X^2).$ Then $R$ is a Gorenstein local ring of dimension one with the maximal ideal $\frak m=(X,Y)R$. Assume that $\frak p=XR$. Then we have the exact sequence of $R$-modules $0\To R/\frak p\To G(R/\frak p)\To E\To 0$  such that $E$ is injective. We notice that $\Gamma_{\frak m}(R/\frak p)=0$ and so we show that $\Gamma_{\frak m}(G(R/\frak p))$ is a  non-zero module. Otherwise, applying the functor $\Gamma_{\frak m}(-)$ to the above exact sequence and using \cref{secgor}, we have $\Gamma_{\frak m}(E)\cong H_{\frak m}^1(R/\frak p)$. Since $\overline{R}=R/\frak p\cong k\llbracket Y\rrbracket$ is a regular local ring of dimension one, we have $H_{\frak m}^1(R/\frak p)\cong E_{\overline{R}}(k)\cong\Hom_R(R/\frak p,E(k))$. Applying $\Hom_R(k,-)$ to these isomorphisms, we deduce that $\Gamma_{\frak m}(E)=E(k)$ so that $E(k)\cong\Hom_R(R/\frak p,E(k))$. But this implies that $\frak p=0$ which is a contradiction.    
\end{Example}

\medskip

\begin{Proposition}\label{evg}
Let $M$ be an $R$-module. Then every injective extension of $M$ can be extended to a Gorenstein injective preenvelop of $M$.
\end{Proposition}
\begin{proof}
Let $M\subset E$ be an injective extsenion of $M$ and so there is an exact sequence of $R$-modules $0\To M\To E\To E/M\To 0$. According to [K, Theorem 7.12], the module $E/M$ fits into an exact sequence of $R$-modules $0\To H\To X\To E/M\To 0$ such that $X\in^{\bot}\cG\cI$ and $H$ is Gorenstein injective. Then we have the commutative diagram 
$$\xymatrix{0\ar[r]&M\ar[r]\ar[d]& E\ar[d]\ar[r]&E/M\ar@{=}[d] \ar[r]& 0\\
0\ar[r]&G\ar[r]&X \ar[r]& E/M\ar[r]& 0.}$$
which gives rise to an exact sequence of modules
 $0\To M\To E\oplus G\To X\To 0.$  
\end{proof}

\medskip

\begin{Proposition}\label{tor}
Let $\frak a$ be an ideal of $R$. If $M$ is an $\frak a$-torsion module, then so is $G(M)$. 
\end{Proposition}
\begin{proof}
By \cref{secgor}, the module $\Gamma_{\frak a}(G(M))$ is Gorenstein injective. Denoted by $i:M\To G(M)$ we show that $\Gamma_{\frak a}(i):M\To \Gamma_{\frak a}(G(M))$ is a Gorenstein injective preenvelope of $M$. Given a Gorenstein injective $R$-module $H$ and a homomorphism $\theta:M\To H$, it factors to $M\stackrel{\Gamma_{\frak a}(\theta)}\To\Gamma_{\frak a}(H)\stackrel{j}\To H$ such that $\Gamma_{\frak a}(j)=1_{\Gamma_{\frak a}(H)}$. Since $G(M)$ is the Gorenstein injective envelope of $M$, there exists a homomorphism $\beta:G(M)\To H$ such that $\beta i=\theta$. Now, considering the homomorphism $j\Gamma_{\frak a}(\beta):\Gamma_{\frak a}(G(M))\To H$, we have $j\Gamma_{\frak a}(\beta)\Gamma_{\frak a}(i)=j\Gamma_{\frak a}(\theta)=\theta.$ Since $\Gamma_{\frak a}(i):M\To\Gamma_{\frak a}(G(M))$ is a Gorenstein injective preenvelope of $M$, the module $G(M)$ is a direct summand of $\Gamma_{\frak a}(G(M))$ which implies that $G(M)=\Gamma_{\frak a}(G(M))$.
\end{proof}


\medskip
\begin{Proposition}\label{go}
Let $G$ be a Gorenstein injective $R$-module. Then $G\cong E\oplus H$ where $E$ is injective and $H$ is reduced Gorenstein injective.
\end{Proposition}
\begin{proof}
 If $G$ is reduced, there is nothing to prove; otherwise the set of all injective submodules $E$ of $G$ is non-empty and since $R$ is noethreian it satisfies the Zorn's lemma conditions under the inclusion. This implies that $G$ has a maximal injective submodule $E$. Setting $H=G/E$, the maximality of $E$ implies that $H$ is reduced. Now, since $R$ has a dualizing complex, by [CFH, Theorem 6.9] the class of Gorenstein injective modules closed under arbitrary direct sums and so by [H1, Proposition 1.4], it is closed under direct summands; consequently $H$ is Gorenstein injective.
\end{proof}

\medskip

\begin{Corollary}
  Let $M$ be an $R$-module and let $M\To G$ be a Gorenstein injective preenvelope of $M$. Then  $G$ has a reduced Gorenstein injective direct summand $H$ such that $H=G(N)$ for some submodule $N$ of $M$ with $M/N$ injective.
\end{Corollary}
\begin{proof}
 If $G$ is reduced there in nothing to prove. If $G$ is not reduced, by \cref{go}, there exists a Gorenstein injective module $H$ and an injective module $E$ such that $G=H\oplus E$. We have the pullback diagram

$$\xymatrix{&0\ar[d]&0\ar[d]\\
0\ar[r]&N\ar[r]\ar[d]& H\ar[d]\ar[r]&X\ar@{=}[d] \ar[r]& 0\\
0\ar[r]&M\ar[d]\ar[r]&G\ar[r]\ar[d]& X\ar[r]& 0\\
&E\ar@{=}[r]\ar[d]&E\ar[d]\\
&0&0}$$
which implies that $N\To H$ is a Gorenstein injective preenvelope of $N$ and since $H$ is reduced, \cref{prred} implies that $H$ is the  Gorenstein envelope of $N$.  
\end{proof}

\section{Localization of indecomposable Gorenstein injective  modules}
We start this section with a result which has already been proved by Enochs and Jenda [EJ1] for Gorenstein rings. We show that it holds without this condition.

\medskip
\begin{Proposition}\label{ind}
Let $\frak p$ be a prime ideal of $R$. Then $G(R/\frak p)$ is indecomposable. Moreover, if $G(R/\frak p)$ is not reduced, then $G(R/\frak p)\cong E(R/\frak p)$ and $R/\frak p\in^{\bot}\cG\cI$.
.
\end{Proposition}
\begin{proof}
Taking  $M=R/\frak p$ and $G=G(M)=G_1\oplus G_2$ and, it follows from [H1, Proposition 1.4] that $G_1$ and $G_2$ are Gorenstein injctive. If $G_1\cap M=0$, then $G_1\oplus M\subset G_1\oplus G_2$ and so $M\subset G_2$ so that $G/M=G_1\oplus (G_2/M)$ and since $G/M\in^{\bot}\cG\cI$, we conclude that $G_2/M\in^{\bot}\cG\cI$. This fact implies that $G_2$ is a Gorenstein injective preenvelope of $M$ so that $G$ is a direct summand of $G_2$; consequently $G_1=0$ which is a contradiction. A similar proof deduces that $G_2\cap M\neq 0$. But $(G_1\cap M)\cap (G_2\cap M)=0$ contradicts the fact that $M$ is domain. If $G$ is not reduced, then it contains an injective submodule $E$ and since $E\in^{\bot}\cG\cI$ and $G$ is a Gorenstein essential extension of $M$, we have $E\cap M\neq 0$ so that $\frak p\in\Ass_RE$. Therefore $E(M)$ is a direct summand of $E$ and so is a direct summand of $G$. Since $G$ is indecomposable, we have $G(M)\cong E(M)$. As $G$ and $G/M$ are in $^{\bot}\cG\cI$, we have  $M\in^{\bot}\cG\cI$.   
\end{proof}

The following examples show if $G$ is an indecomposable Gorenstein injective module, then it may not exist a prime ideal $\frak p\in\Spec R$ such that  $G=G(R/\frak p)$. 

\medskip
\begin{Example}\label{ex}
Let $(R,\frak m)$ be a local ring and $x,y\in\frak m$ such that $xR=\Ann_R(y)$ and $yR=\Ann_R(x)$ and let $a\in\frak m$ be a regular element on $R/(x,y)$. Set 
$\gamma_a=\begin{pmatrix}x&a\\0&y\end{pmatrix}$ and
$\eta_a=\begin{pmatrix}y&-a\\0&x\end{pmatrix}$. By [H2, Lemma 3.4], there is an exact complex of free $R$-modules 
$$\cF:\dots \To R^2\stackrel{\gamma_a}\To R^2\stackrel{\eta_a}\To R^2\stackrel{\gamma_a}\To R^2\stackrel{\eta_a}\To\dots $$ such that $\Hom_R(\cF,A)$ is exact. Hence $G_a=\Coker\gamma_a$ and $H_a=\Coker \eta_a$ are finitely generated Gorenstein projective. It follows from [H2, Corollary 3.9, Theorem 4.5, Theorem 5.1] that $G_a$ and $H_a$ are non-free, non-isomorphic and indecomposable. Setting $G(a)=\Hom_R(G_a,E(R/\frak m))$ and $H(a)=\Hom_R(H_a,E(R/\frak m))$, they are artinian and Gorenstein injective by [CFH, Proposition 5.1]. Moreover, $\End_{R}(G(a))\cong \widehat{R}\cong \End_R(H(a))$ is a local ring so that $G(a)$ and $H(\frak a)$ are indecomposable. Further, the faithfully flatness of $\widehat{R}$ implies that $G(a)$ and $H(a)$ are not injective and they are non-isomorphic. Consequently,  at least one of them is not $G(R/\frak m)$.  
\end{Example}

 \medskip

\begin{Example}\label{eexx}
Let $k$ be a field and set $R=k\llbracket X,Y\rrbracket/(XY).$ Then $R$ is a Gorenstein local ring of dimension one with the maximal ideal $\frak m=(X,Y)R$. Then $R/XR$ and $R/YR$ are finitely generated Gorenstein projective and so $G=\Hom_R(R/XR,E(k))$ and $H=\Hom_R(R/YR,E(k))$ are artinian Gorenstein injective $R$-modules and since $\End_R(G)\cong R/XR$ and $\End_R(H)\cong R/YR$ are local rings, they are  indecomposable. We now prove the following assertions. 

{\bf Claim 1}. $G(k)\neq G,H $. We prove the claim for $G$ and for $H$ is similarly. Otherwise, since $R$ is a Gorenstein ring of dimension one, there exists an exact sequence $0\To k\To G\To I\To 0$ of $R$-modules such that $I$ is injective. According to [FS, Example 5.2], there exists an exact sequence of modules $0\To k\To G\stackrel{Y.}\To G\To 0$ which yields a commutative diagram with exact rows  
$$\xymatrix{0\ar[r]&k\ar[r]\ar@{=}[d]& G\ar[d]\ar[r]&I\ar[d] \ar[r]& 0\\
0\ar[r]&k\ar[r]&G \ar[r]& G\ar[r]& 0.}$$
 and so an exact sequence of modules $0\To G\To I\oplus G\To G\To 0$. We observe that $\Hom_R(k,G)\cong k$ and so applying $\Hom_R(k,-)$, we have an exact sequence of $k$-vector spaces $0\To k\To \Hom_R(k,I)\oplus k\To k$ which implies that $\Hom_R(k,I)\cong k$  so that $I\cong E(k)$. Now, the following commutative diagram with exact rows
 
 $$\xymatrix{0\ar[r]&k\ar[r]\ar[d]^0& G\ar[d]^{Y.}\ar[r]&E(k)\ar[d]^{Y.} \ar[r]& 0\\
0\ar[r]&k\ar[r]&G \ar[r]& E(k)\ar[r]& 0.}$$ 
 which forces an exact sequence of modules  $0\To k\To E(k)\To E(k)\To 0$. But this implies that $R$ is regular; which is a contradiction.

{\bf Claim 2}. $G(k)$ is a direct summand of $H\oplus G$. It is clear that $$\xymatrix{\dots\ar[rr]^{\begin{pmatrix}Y&0\\0&X\end{pmatrix}}&&R^2
\ar[rr]^{\begin{pmatrix}X&0\\0&Y\end{pmatrix}} &&R^2 \ar[rr]^{\begin{pmatrix}Y&0\\0&X\end{pmatrix}}\ar[rr] &&R^2 \ar[rr]^{\begin{pmatrix}X&Y\end{pmatrix}} &&R\ar[r]&k\ar[r]& 0}$$
is a free resolution of $k$ and $L=\Coker\begin{pmatrix}Y&0\\0&X\end{pmatrix}=R/YR\oplus R/XR$. Hence applying $\Hom_R(-,E(k)$, we have an exact sequence $0\To k\To E(k)\To G\oplus H\To 0$. On the other hand, if we apply $\Hom_R(-,E(k))$ to the exact sequence $\dots\stackrel{X.}\To R\stackrel{Y.}\To R\stackrel{X.}\To R\stackrel{Y.}\To \dots$, we have  two exact sequence of $R$-modules $0\To G\To E(k)\To H\To 0$ and $0\To H\To E(k)\To G\To 0$ and so an exact sequence $0\To H\oplus G\To E(k)\oplus E(k)\To G\oplus H\To 0$. Thus the following commutative diagram with exact rows 
 $$\xymatrix{0\ar[r]&k\ar[r]\ar[d]& E(k)\ar[d]\ar[r]&G\oplus H\ar@{=}[d] \ar[r]& 0\\
0\ar[r]&H\oplus G\ar[r]&E(k)\oplus E(k) \ar[r]&  G\oplus H\ar[r]& 0}$$ 
 gives rises to an exact sequence of modules $0\To k\To H\oplus G\To E(k)\To  0$ which implies that $k\To H\oplus G$ is a Gorenstein injective preenvelope of $k$; and consequently $G(k)$ is a direct summand of $H\oplus G$.

{\bf Claim 3}. $\mu(\frak m,G(k))=1$. It follows from \cref{prred} and Claim 2 that $G(k)\oplus E\cong H\oplus G$ for some injective module $E$.
Applying $\Hom_R(k,-)$ and the fact that $G(k)$ is indecomposable, we deduce that $\Hom_R(k,G(k))\cong k$.

\end{Example}

\medskip

\begin{Lemma}\label{gorr}
Let $\frak p$ be a prime ideal of $R$ and $G$ be a Gorenstein injective $R_{\frak p}$-module. Then $G$ is a Gorenstein injective $R$-module.
\end{Lemma}
\begin{proof}
Given a complete $R_{\frak p}$-injective resolution $\cI:\dots \To I_i\To I_{i-1}\To\dots $ of $G$, since for each $i$, we have $\Hom_R(I_i,I_{i-1})=\Hom_{R_{\frak p}}(I_i,I_{i-1})$, $\cI$ is an exact sequence of injective $R$-modules. For arbitrary injective $R$-module $E$, there is a decomposition $E=E'\oplus E''$ for $E$ such that $E'=\bigoplus_{\frak q\subset \frak p}E(R/\frak q)$ and $E''=\bigoplus_{\frak q\nsubseteq \frak p}E(R/\frak q)$. Then for each $i$, we have $\Hom_R(E'',I_i)=\prod_{\frak q\nsubseteq \frak p}\Hom_R(E(R/\frak q), I_i)$. Fixing a prime ideal $\frak q$ such that  $\frak q\nsubseteq \frak p$ and an element $s\in \frak q\setminus \frak p$, the map $I_i\stackrel{s.}\To I_i$ is isomorphism and so is $\Hom_R(E(R/\frak q),I_i)\stackrel{s.}\To \Hom_R(E(R/\frak q),I_i)$. But for every $f\in\Hom_R(E(R/\frak q),I_i)$ and $x\in E(R/\frak q)$, the element $f(x)$ is annihilated by some power of $s$ and so $f(x)=0$. Thus $\Hom_R(E(R/\frak q),I_i)=0$ and consequently $\Hom_R(E'',I_i)=0$. Therefore $\Hom_R(E,\cI)=\Hom_R(E',\cI)=\Hom_{R_{\frak p}}(E'_{\frak p},\cI)$ is exact. This implies that $G$ is a Gorenstein injective $R$-module.  
\end{proof}

\medskip
\begin{Lemma}\label{preen}
Let $\frak p\in\Spec R$ and let $M$ be an $R$-module with the Gorenstein envelope $\iota:M\To G(M)$. Then $\iota_{\frak p}$ is a Gorenstein injective  preenvelope of $R_{\frak p}$-module $M_{\frak p}$.
\end{Lemma}
\begin{proof}
There is an exact sequence of $R$-modules $0\To M\stackrel{\iota}\To G(M)\To X\To 0$ such that $X\in^{\bot}\cG\cI$. For any Gorenstein injective $R_{\frak p}$-module $H$, assume that $\cI:\dots\To I_i\To I_{i-1}\To \dots$ is a complete $R_{\frak p}$-injective resolution of $H$. Then there is an isomorphism $\Hom_{R_{\frak p}}(X_{\frak p}, \cI)\cong\Hom_R(X,\cI)$ and since by \cref{gorr}, the $R$-module $H$ is Gorenstein injective, $\Hom_{R_{\frak p}}(X_{\frak p}, \cI)$ is exact so that $\Ext_{R_{\frak p}}^1(X_{\frak p},H)=0$. Thus $\iota_{\frak p}:M_{\frak p}\To G(M)_{\frak p}$ is a Gorenstein injective preenvelope of $R_{\frak p}$-module $M_{\frak p}.$
\end{proof}

For every $\frak p \in\Spec R$, set $V(\frak p)=\{\frak q\in\Spec R|\hspace{0.1cm} \frak p\subset \frak q\}$ and for every integer $k$, we define the subset $X_k(\frak p)$ of $\Spec R$ as $X_k(\frak p)=\{\frak q\in V(\frak p)|\hspace{0.1cm} {\rm ht}\frak p-\sup D_{\frak p}=k\}$. For every $\frak q\in \Spec R$, we have the following (in)equalities  $$\id_D\geq \id_{R_{\frak q}}D_{\frak q}=\id_{R_{\frak q}}\Sigma^{{\rm ht}\frak q-\sup D_{\frak q}}D_{\frak q}+{\rm ht}\frak q-\sup D_{\frak q}={\rm ht} {\frak q}-\sup D_{\frak q}\geq {\rm ht}\frak q-\sup D\geq -\sup D$$
where the second equality holds as $\Sigma^{{\rm ht}\frak q-\sup D_{\frak q}}D_{\frak q}$  is a normalized dualizing complex. We find out from this description that there are only finitely many such $X_k$. As $\Sigma^nD$ is a dualizing complex for every integer $n$, we may assume that $\sup D=0$ and  we have the following theorem.

\medskip
\begin{Theorem}\label{filt}
Let $\frak p\in\Spec R$ such that $\height \frak p-\sup D_{\frak p}=t$ and $\sup D=0$. Then there exists a finite filtration $0\subset G_d\subset G_{d-1}\subset\dots \subset G_t=G(R/\frak p)$ of Gorenstein injective submodules of $G(R/\frak p)$ such that $G_k/G_{k+1}\cong\bigoplus_{\frak q\in X_{k}(\frak p)}\Tor_{k}^R(E(R/\frak q),\RHom_R(D,G(R/\frak p))$ is reduced for each $t\leq k\leq d$. In particular, 
$$G_t/G_{t+1}\cong \Tor_t^R(E(R/\frak p),\RHom_R(D,G(R/\frak p)))\cong G(R/\frak p)_{\frak p}.$$
\end{Theorem}
\begin{proof}
Without loss of generality, we may assume that $\sup D=0$ and assume that $\id_RD=d$. By [S2, Theorem 2.12], the module $G(R/\frak p)$ admits a finite filtration of Gorenstein injective submodules $$0=G_{d+1}\subset G_{d}\subset\dots \subset G_{1}\subset G_{0}=G(R/\frak p)$$ such that each $G_{k}/G_{k+1}\cong \bigoplus_{\frak q\in X_{k}}\Tor_{k}^R(E(R/\frak q),\RHom_R(D,G(R/\frak p)))$ is Gorenstein injective. Fix $0\leq k\leq d$ and $\frak q\in X_k$. If $\frak p\nsubseteq\frak q$, there exists $x\in\frak p\setminus \frak q$ and so the isomorphism $E(R/\frak q)\stackrel{x.}\To E(R/\frak q)$ give rises to an isomorphism  $$\Tor_{k}^R(E(R/\frak q),\RHom_R(D,G(R/\frak p)))\stackrel{x.}\To \Tor_{k}^R(E(R/\frak q),\RHom_R(D,G(R/\frak p))).$$ But since 
$G(R/\frak p)$ is $\frak p$-torsion, we deduce that 
$\Tor_{k}^R(E(R/\frak q),\RHom_R(D,G(R/\frak p)))=0$. Thus $G_{k}/G_{k+1}\cong \bigoplus_{\frak q\in X_{k}(\frak p)}\Tor_{k}^R(E(R/\frak q),\RHom_R(D,G(R/\frak p))).$ On the other hand, for each $k<t$ and $\frak  q\in X_k(\frak p)$, we have $\mu_{R_{\frak q}}^i(D_{\frak q})\neq 0$ for all $i\neq k$ and so it follows from [CF, Lemma 6.1.19] that $\mu_{R_{\frak p}}^{i-\dim R_{\frak q}/\frak pR_{\frak q}}(D_{\frak p})=0$ for all $i\neq k$ which implies that $t=k-\dim R_{\frak q}/\frak pR_{\frak q}=k+{\rm ht}\frak p-{\rm ht}\frak q<k$ which is a contradiction as $k<t$; the first equality holds as $\mu_{R_{\frak p}}^t(D_{\frak p})=1$ and $\mu_{R_{\frak p}}^i(D_{\frak p})=0$ for all $i\neq t$ and the second equality holds as $R_{\frak q}$ is a catenary ring by [CF, Proposition 7.1.12]. Therefore $X_k(\frak p)$ is an empty set and so $G_k/G_{k+1}=0$ for each $k<t$ so that $G(R/\frak p)=G_0=\dots=G_t$. To do the second assertion, if $\frak q\in X_t(\frak p)$ such that $\frak p\subsetneq \frak q$, then $\mu_{R_{\frak q}}^i(D_{\frak q})=0$ for all $i\neq t$ and so by a similar argument as mentioned before, we have $\mu_{R_{\frak p}}^j(D_{\frak p})= 0$ for all $j\neq t-{\rm ht}\frak p+{\rm ht}\frak q$ which is a contradiction. Moreover, by [S2, Lemma 2.6 and Proposition 2.5] and the fact that $G(R/\frak p)$ is $\frak p$-torsion, we deduce that $G_t/G_{t+1}\cong G(R/\frak p)_{\frak p}$. Fixing $t\leq k\leq d$, we prove that $G_{k}/G_{k+1}$ is reduced. Otherwise, $X=G_{k}/G_{k+1}$ has an injective submodule $E$ and so we have the following pull back diagram 

$$\xymatrix{&&0\ar[d]&0\ar[d]\\
0\ar[r]&G_{k+1}\ar@{=}[d]\ar[r]& D\ar[d]\ar[r]&E \ar[r]\ar[d]& 0\\
0\ar[r]&G_{k+1}\ar[r]&G_{k}\ar[r]\ar[d]& X\ar[r]\ar[d]& 0\\
&&Y\ar@{=}[r]\ar[d]&Y\ar[d]\\
&&0&0}$$
which implies that the top row splits. But, since $G(R/\frak p)$ is indecomposable, $G_k$ is reduced and so is $D$; which is a contradiction.  
 \end{proof}
\medskip

\begin{Corollary}\label{corlocc}
Let $\frak p\in\Spec R$. Then $G(R/\frak p)_{\frak p}\cong G(k(\frak p))$, where $k(\frak p)=R_{\frak p}/\frak pR_{\frak p}$.
\end{Corollary}
\begin{proof}
Since $\iota:R/\frak p\To G(R/\frak p)$ is the Gorenstein envelope of $G(R/\frak p)$, the map $\iota_{\frak p}$ is a Gorenstein injective preenvelope of $k(\frak p)$ and so using \cref{prred} there exists an injective $R_{\frak p}$-module $E$ such that $G(R/\frak p)_{\frak p}=G(R_{\frak p}/\frak pR_{\frak p})\oplus E$ and since by \cref{filt} the Gorenstein injective module $G(R/\frak p)_{\frak p}$ is reduced, we have $E=0$. 
\end{proof}
\medskip

Suppose that ${\bf x}=x_1,\dots, x_t$ is a sequence of elements of $R$. The Koszul complex of $R$ induced by ${\bf x}$ is denoted by $K_{\bullet}({\bf x})$. For any $R$-module $M$, we denote the complex $\Hom_A(K_{\bullet}({\bf x}),M)$ by 
$K^{\bullet}({\bf x},M)$ and for each $i\geq 0$, we denote $H^i(K^{\bullet}({\bf x},M))$ by $H^i({\bf x},M)$

\begin{Proposition}\label{kos}
Let $G$ be a Gorenstein injective $R$-module and $\frak a=({\bf x})$ be an ideal of $R$ such that $H^0({\bf x},G)=0$. Then $H^i({\bf x},G)=0$ for all $i>0$.
\end{Proposition}
 \begin{proof}Since $(0:_G\frak a)=H^0({\bf x},G)=0$, applying the functor $\Hom_R(-,G)$ to the exact sequences $0\To \frak a/\frak a^2\To R/\frak a^2\To R/\frak a\To 0$ and  $0\To L_1\To (R/\frak a)^{t_1}\To \frak a^2/\frak a\To 0$, we deduce that $(0:_G\frak a^2)=0$. Continuing this way for the  exact sequences of $R$-modules $0\To \frak a^n/\frak a^{n-1}\To R/\frak a^n\To R/\frak a^{n-1}\To 0$ and $0\To L_{n-1}\To (R/\frak a)^{t_{n-1}}\To \frak a^n/\frak a^{n-1}\To 0$ , we conclude that $(0:_G\frak a^n)=0$ for all $n\geq 1$; and consequently $\Gamma_{\frak a}(G)=0$. There exists an exact sequence of modules $0\To G\To E\To G_1\To 0$ such that $E$ is the injective envelope of $G$ and $G_1$ is Gorenstein injective. Since $\Gamma_{\frak a}(G)=0$, we deduce that $\Gamma_{\frak a}(E)=0$ and so the previous exact sequence and \cref{secgor} imply that $\Gamma_{\frak a}(G_1)=0$.    
 The preceding  argument implies that 
 $\Hom_R(R/\frak a,E)=0=\Hom_R(R/\frak a,G_1)$. 
 It now follows from [AM, Lemma 2.8] that $H^i({\bf x},E)=0$ for all $i>0$. Moreover, the exact sequence of complexes $0\To K({\bf x},G)\To K({\bf x},E)\To K({\bf x},G_1)\To 0$ implies that $H^{i+1}({\bf x},G)\cong H^i({\bf x},G_1)$ for all $i\geq 0$. Now an easy induction on $i$ implies that $H^i({\bf x},G)=0$ for all $i\geq 0$. 
\end{proof}

\medskip
\begin{Corollary}\label{loc}
Let $\frak p\in\Spec R$ and let $G$ be a Gorenstein injective $R$-module such that $\Ass_RG=\{\frak p\}$. Then the map $G\stackrel{s.}\to G$ is isomorphism for all $s\in R\setminus \frak p.$ Furthermore, $G$ is an $R_{\frak p}$-module.              
\end{Corollary}
\begin{proof}
Fix $s\in R\setminus \frak p.$ Since $\Ass_RG=\{\frak p\}$, we have $E(G)=E(R/\frak p)^{(X)}$ for some set $X$ and hence the map $E(G)\stackrel{s.}\to E(G)$ is isomorphism so that $G\stackrel{s.}\to G$ is injective.  Thus $H^0(s,G)=0$; and consequently $H^1(s,G)=\Coker(s.)=0$ by \cref{kos}. To prove the second claim for each $r/s\in R_{\frak p}$ and $g\in G$, by the first case there exists a unique element $g_1\in G$ such that $g=sg_1$ and so we define $(r/s)g:=rg_1$.
\end{proof}

 \medskip
\begin{Theorem}\label{locgor}
Let $\frak p\in\Spec R$ with ${\rm ht}\frak p-\sup D_{\frak p}=t$ and $\sup D=0$. Then the following conditions are equivalent.

${\rm (i)}$ $\Ass_R G(R/\frak p)=\{\frak p\}$.

${\rm (ii)}$ $G(R/\frak p)$ is an $R_{\frak p}$-module.

${\rm (iii)}$ $G(R/\frak p)=G(k(\frak p))$ where $k(\frak p)=R_{\frak p}/\frak pR_{\frak p}$.

${\rm (iv)}$ $\Tor_i^R(E,\RHom_R(D,G(R/\frak p)))=0$ for all $i>t$ and all injective modules $E$ .
\end{Theorem}
\begin{proof}
(i)$\Rightarrow$ (ii) is \cref{loc}. (ii)$\Rightarrow$ (iii) is obtained from  \cref{corlocc}. (iii)$\Rightarrow$ (i). $G(R/\frak p)$ is embedded in $E(R/\frak p)^{(X)}$ and hence $\Ass_RG(R/\frak p)=\{\frak p\}$. (i)$\Rightarrow$ (iv). As $E=\bigoplus_{\frak q\in\Spec R}E(R/\frak q)^{\mu_{\frak q}}$, it suffices to consider $E=E(R/\frak q)$ and we show that $\Tor_i^R(E(R/\frak q),\RHom_R(D,G(R/\frak p)))=0$ for all $i>t$. Fixing $i>t$, if $\frak p\nsubseteq \frak q$, then there is $x\in\frak p\setminus\frak q$ so that $E(R/\frak q)\stackrel{x.}\To E(R/\frak q)$ is isomorphism which gives rises to an isomorphism  $\Tor_i^R(E(R/\frak q),\RHom_R(D,G(R/\frak p)))\stackrel{x.}\To\Tor_i^R(E(R/\frak q),\RHom_R(D,G(R/\frak p)))$. Thus  $\Tor_i^R(E(R/\frak q),\RHom_R(D,G(R/\frak p)))=0$ as $G(R/\frak p)$ is $\frak p$-torsion. If $\frak p\subsetneq \frak q$, then there is $x\in \frak q\setminus \frak p$ and so by \cref{loc}, the map $G(R/\frak p)\stackrel{x.} \To G(R/\frak p)$ is isomorphism and since $E(R/\frak q)$ is $\frak q$-torsion, a similar argument mentioned  before implies that $\Tor_i^R(E(R/\frak q),\RHom_R(D,G(R/\frak p)))=0$. The case $\frak q=\frak p$ follows from [S2, Lemma 2.4]. (iv)$\Rightarrow$(ii). It follows from \cref{filt} that $G_d=\dots=G_{t+1}=0$ and so $G(R/\frak p)=G_t=\Tor_t^R(E(R/\frak p),\RHom(D,G))\cong\Gamma_{\frak pR_{\frak p}}(G(R/\frak p)_{\frak p})\cong G(R/\frak p)_{\frak p}$; where the first isomorphism holds by  [S2, Proposition 2.5] and the second isomorphism holds because $G(R/\frak p)$ is $\frak p$-torsion.
\end{proof}

 The following example shows that for every $\frak p\in\Spec R$, the Grenstein injective module $G(R/\frak p)$ is not an $R_{\frak p}$-module in general even if $R$ is a Gorenstein local ring of dimension one. 
 
\medskip
  \begin{Example}
  Let $k$ be a field and $R=k\llbracket X,Y\rrbracket/(X^2).$ Then $R$ is a Gorenstein local ring of dimension one with the maximal ideal $\frak m=(X,Y)R$ and set $\frak p=XR$. There exists an exact sequence of $R$-modules $\To R/\frak p\To G(R/\frak p)\To E\To 0$ such that $E$ is injective and by \cref{exx}, we have $\Gamma_{\frak m}(G(R/\frak p))$ is non-zero which implies that $\frak m\in\Ass_RG(R/\frak p)$. It now follows from \cref{locgor} that $G(R/\frak p)$ is not an $R_{\frak p}$-module.  
    
  \end{Example}

\section{Bass numbers of indecomposable Gorenstein injective modules}
In this section we study the Bass numbers of $G(R/\frak p)$ for a prime ideal $\frak p\in\Spec R$. We start with the following lemma.

\begin{Lemma}\label{injcover}
Let $\frak a$ be an ideal of $R$ and let $G$ be an $\frak a$-torsion Gorenstein injective $R$-module. Then the injective cover of $G$ is $\frak a$-torsion.
\end{Lemma}
\begin{proof}
Assume that $E\To G$ is the injective cover of $G$. Then it fits into an exact sequence of $R$-modules $0\To K\To E\To G\To 0$ such that $K$ is Gorenstein injective. Application of $\Gamma_{\frak a}(-)$ to this exact sequence induces the exact sequence of modules $0\To \Gamma_{\frak a}(K)\To \Gamma_{\frak a}(E)\To G\To 0$. By \cref{secgor}, the module $\Gamma_{\frak a}(K)$ is Gorenstein injective so that $\Gamma_{\frak a}(E)\To G$ is an injective precover of $G$. Therefore $E$ is a direct summand of $\Gamma_{\frak a}(E)$ which implies that $\Gamma_{\frak a}(E)=E$.
\end{proof} 

\begin{Corollary}
Let $\frak p\in\Spec R$. Then the injective cover of $G(R/\frak p)$ is $\frak p$-torsion. 
\end{Corollary}
\begin{proof}
$G(R/\frak p)$ is $\frak p$-torsion by \cref{tor} and hence the result follows by \cref{injcover}.
\end{proof}

\medskip
\begin{Proposition}\label{gorred}
Let $G$ be a Gorenstein injective $R$-module. Then $E(G)/G$ is reduced. Moreover, If $G$ is reduced, then $E(G)\To E(G)/G$ is the injective cover of $E(G)/G$.  
\end{Proposition}
\begin{proof}
Set $H=E(G)/G$ and and let $E$ be a non-zero injective submodule of $H$. Then we have the following commutative diagram
$$\xymatrix{&&&E\ar@{^{(}->}[ld]^g\ar@{^{(}->}[d]^j\\
0\ar[r]&G\ar[r]^i&E(G)\ar@{->>}[ld]^p \ar[r]^f& H\ar[r]& 0\\&E_1}$$
where the existense of $g$ is obtained from the fact that $f$ is an injective precover and it is monic as $j$ is monic. We also notice that $E\hookrightarrow E(G)\twoheadrightarrow E_1$ is a split short exact sequence and hence the diagram forces that $pi$ is monic which contradicts the fact that $E(G)$ is the injective envelope of $G$. Assume that $E\To H$ is the injective cover of $H$ and so we have the following commutative diagram $$\xymatrix{0\ar[r]&M\ar[r]\ar@{^{(}->}[d]& E\ar@{^{(}->}[d]\ar[r]&H\ar@{=}[d] \ar[r]& 0\\
0\ar[r]&G\ar[r]&E(G) \ar[r]& H\ar[r]& 0.}$$ 
 Since $E(G)\To H$ is an injective precover of $H$, the injective module $E$ is a direct summand of $E(G)$ and so $M$ is a direct summand of $G$ as well.  Consequently, the diagram forces that $E(G)/E\cong G/M$ is an injective submodule of $G$ and since $G$ is reduced, we have $E=E(G)$. 
\end{proof}
\medskip

\begin{Lemma}\label{corlocc}
Let $\frak p\in\Spec R$. Then $\mu^i(\frak q, G(R/\frak p))=\mu^i(\frak qR_{\frak q},G(R_{\frak q}/\frak pR_{\frak q}))$ for all $\frak q\in\Spec R$ and all $i\geq 0$.
\end{Lemma}
\begin{proof}
If $\frak p\nsubseteq \frak q$, then $ G(R/\frak p)_{\frak q}=0$ as $G(R/\frak p)$ is $\frak p$-torsion and so there is nothing to prove in this case. If $\frak p\subset \frak q$, by \cref{preen}, since $\iota:R/\frak p\To G(R/\frak p)$ is the Gorenstein envelope of $G(R/\frak p)$, the map $\iota_{\frak q}$ is a Gorenstein injective preenvelope of $R_{\frak q}/\frak pR_{\frak q}$ and so using \cref{prred} there exists an injective $R_{\frak q}$-module $E$ such that $G(R/\frak p)_{\frak q}=G(R_{\frak q}/\frak pR_{\frak q})\oplus E$. On the other hand, it follows from \cref{filt} that there exists an integer $k$ such that $\frak q\in X_k(\frak p)$ and so by [S2, Lemma 2.6 and Proposition 2.5], we have $\Gamma_{\frak qR_{\frak q}}(G(R/\frak p)_{\frak q})\cong \Tor_{k}^R(E(R/\frak q),\RHom_R(D,G(R/\frak p))$ which is reduced by \cref{filt}. Therefore $\Gamma_{\frak q}(E)=0$; and consequently $\Ext^i_{R_{\frak q}}(k(\frak q),G(R/\frak p)_{\frak q})\cong \Ext^i_{R_{\frak q}}(k(\frak q),G(R_{\frak q}/\frak pR_{\frak q}))$ for each $i\geq 0$. 
\end{proof}
\medskip

\medskip
\begin{Theorem}\label{gid}
Let $\frak p\in\Spec R$ such that $\Gid_RR/\frak p$ is finite. Then $\mu^i(\frak q,G(R/\frak p))$ is finite for all $i\geq 0$ and all $\frak q\in\Spec R$.
\end{Theorem}
\begin{proof}
Given $\frak q\in\Spec R$, by \cref{corlocc}, we have  $\mu^i(\frak q,G(R/\frak p))=\mu^i(\frak qR_{\frak q},G(R_{\frak q}/\frak pR_{\frak q}))$ for all $i\geq 0$ and since by \cref{tor}, the module $G(R/\frak p)$ is $\frak p$-torsion, we may assume that $\frak p\subset \frak q$ . Then replacing $R$ by the local ring $(R_{\frak q},\frak qR_{\frak q})$, we may assume that $R$ is a local ring with the maximal ideal $\frak m=\frak q$. Taking $\Gid_RR/\frak p=n$ and considering the minimal injective resolution $0\To R/\frak p\To E^0\To E^1\To\dots$ of $R/\frak p$, for each $i\geq n$, the module $G^i=\Ker(E^i\To E^{i+1})$ is Gorenstein injective. It follows from \cref{gorred} that $G^{i}$ is reduced for all $i\geq n+1$ and set $G=G^{n+2}$. If $E_0\To G$ is the injective cover of $G$ with $K_1=\Ker(E_0\To G)$, by [EJ2, Proposition 10.1.6], the Gorenstein injective module $K_1$ is reduced and $E_0=E(K_1)$ (we can set $K_0=G$). Replacing $G$ by $K_1$ and  the same argument we have the exact sequence of modules $\dots\To E_1\To E_0\To G\To 0$ such that $K_{i+1}=\Ker(E_{i}\To K_i)$ is reduced Gorenstein injective and $E_i=E(K_{i+1})$ for each $i\geq 0$.
Fixing $i\geq 1$, for the minimal injective resolution  $0\To K_{i}\To E_{i-1}\To E_{i-2}\To\dots$ of $K_{i}$, it follows from \cref{secgor} that $0\To \Gamma_{\frak m}({K_{i}})\To \Gamma_{\frak m}({E_{i-1}})\To \Gamma_{\frak m}({E_{i-2}})\To\dots$ is the minimal injective resolution of $\Gamma_{\frak m}({K_{i}})$. Hence $\Gamma_{\frak m}({K_i})$ is reduced and $\Gamma_{\frak m}({E_i})\To \Gamma_{\frak m}({K_i})$ is the injcetive cover of $\Gamma_{\frak m}({K_i})$ for all $i\geq 0$ by \cref{gorred}. As $R/\frak p$ is finitely generated, $\mu(\frak m,G)$ is finite so that $\Gamma_{\frak m}(G)$ is an artinian $\widehat{R}$-module, where $\widehat{R}$ is the completion of $R$ with respect to $\frak m$. Thus [CFH, Proposition 5.1] implies that $\Hom_R(\Gamma_{\frak m}(G),E(R/\frak m))$ is a finitely generated Gorenstein projective $\widehat{R}$-module so that there exists a positive integer $m$ and an exact sequence of $\widehat{R}$-modules $0\To \Hom_R(\Gamma_{\frak m}(G),E(R/\frak m))\to \widehat{R}^m\To L\To 0$ such that $L$ is finitely generated Gorenstein projective. Application of $\Hom_R(-,E(R/\frak m))$ gives rises to an exact sequence of $\widehat{R}$-modules $0\To\Hom_R(L,E(R/\frak m))\To E(R/\frak m)^m\To \Gamma_{\frak m}(G)\To 0$ such that $\Hom_R(L,E(R/\frak m))$ is an artinian Gorenstein injective $\widehat{R}$-module. It follows from [S1, Lemma 3.5] that $\Hom_R(L,E(R/\frak m))$ is a Gorenstein injective $R$-module and so $E(R/\frak m)^m\To \Gamma_{\frak m}(G)$ is an injective precover of $\Gamma_{\frak m}(G)$. Thus $\Gamma_{\frak m}({E_0})$ is a direct summand of $E(R/\frak m)^m$ because $\Gamma_{\frak m}({E_0})\To\Gamma_{\frak m}(G)$ is the injective cover of $\Gamma_{\frak m}(G)$; consequently only finitely many copies of $E(R/\frak p)$ occurs in $E_0$. Now considering the exact sequences of $R$-mdules $0\To \Gamma_{\frak m}({K_i})\To \Gamma_{\frak m}({E_{i-1}})\To \Gamma_{\frak m}({K_{i-1}})\To 0$  and repeating this argument inductively, we deduce that  $\Gamma_{\frak m}({K_i})$ is artinian for each $i\geq 0$ and so only finitely many copies of $E(R/\frak m)$ occurs in each $E_i$. On the other hand, we have the following commutative diagram of $R$-modules

$$\xymatrix{\cX\ar[d]^f:=&0\ar[r]&R/\frak p\ar[r]\ar[d]^{\epsilon}&E^0\ar[r]\ar[d]^{f_0}&\dots\ar[r]& E^n\ar[r]\ar[d]^{f_n}&G\ar@{=}[d] \ar[r]& 0\\
\cY:=&0\ar[r]&K_{n+1}\ar[r]&E_n\ar[r]&\dots\ar[r]& E_0\ar[r]& G\ar[r]&0.}$$ 
Using con$(f)$, there exists the exact sequences of $R$-modules $0\To R/\frak p\To E^0\oplus K_{n+1}\To X\To 0$ and  $0\To X\To E^1\oplus E_n\To\dots \To E_0\To 0$; consequently $R/\frak p\To E^0\oplus K_{n+1}$ is a Gorenstein injecrive preenvelope of $R/\frak p$. Since $R/\frak p$ is finitely generated, $\mu^i(\frak m,R/\frak p)$ is finite for all $i\geq 0$ and so by the previous description $\mu^i(\frak m,E^0\oplus K_{n+1})$ is finite for all $i\geq 0$ and since $G(R/\frak p)$ is a direct summand of $E^0\oplus K_{n+1}$, we deduce that $\mu^i(\frak m,G(R/\frak p))$ is finite for all $i\geq 0$.   
\end{proof}

The following corollary  is obtained immediately.
\begin{Corollary}\label{goro}
Let $R$ be a Gorenstein ring of Krull dimension $d$. Then $\mu^i(\frak q,G(R/\frak p))$ is finite for all $i\geq 0$ and all $\frak p,\frak q\in\Spec R$. 
\end{Corollary}
\begin{proof}
Since $R$ is Gorenstein, $\Gid_RR/\frak p\leq d$ for all $\frak p\in\Spec R$, and hence the result follows from \cref{gid}. 
\end{proof}

\begin{Corollary}\label{subb}
Let $\frak p\in\Spec R$ such that $\Gid_RR/\frak p$ is finite. Then $\mu^i(\frak r, G(R/\frak q))$  is finite for all $i\geq 0$ and all 
 $\frak r,\frak q\in \Spec R$ with $\frak q\subset\frak p.$
\end{Corollary}
\begin{proof}
Since $\Gid_RR/\frak p$ is finite, by [CFH, Proposition 5.5], we have $\Gid_{R_{\frak p}}R_{\frak p}/\frak pR_{\frak p}$ is finite and so $R_{\frak p}$ is a Gorenstein local ring by [FF, Theorem 4.5]. Now, for every $\frak q\in\Spec R$ with $\frak q\subset \frak p$, the local ring $R_{\frak q}$ is Gorenstein and so for every $i\geq 0$, we have $\mu^i(\frak r, G(R/\frak q))=\mu^i(\frak rR_{\frak q}, G(k(\frak q))$ which is finite by \cref{goro}.
\end{proof}

  \begin{Example}
  Let $k$ be a field and $R=k\llbracket X,Y\rrbracket/(X^2).$ Then $R$ is a Gorenstein local ring of dimension one with the maximal ideal $\frak m=(X,Y)R$ and set $\frak p=XR$. We claim that $\mu^i(\frak p,G(R/\frak p))=1$ for all $i\geq 0$. Consider the exact sequence  of $R$-modules $0\To R/\frak p\To E(R/\frak p)\To G\To 0$. Since $R$ is Gorenstein, $G$ is Gorenstein injective. Since $R_{\frak p}$ is a Gorenstein local ring of dimension zero, $R_\frak p$ is injective and $R_{\frak p}\cong E(k(\frak p))$, where $k(\frak p)=R_{\frak p}/\frak pR_{\frak p}$. Furthermore, the fact that $\frak p^2=0$ implies that $\frak pR_{\frak p}$ is a vector space on $k(\frak p)$ and since $\frak p$ is principal, we have $\frak pR_{\frak p}\cong k(\frak p)$ so that $l(R_{\frak p})=2$. Therefore localizing the above exact sequence at $\frak p$ and using the previous argument, we deduce that $l(G_p)=1$, and hence $G_{\frak p}\cong k(\frak p)$. On the other hand, we have the following commutative diagram of $R$-modules with exact rows
  $$\xymatrix{0\ar[r]&R/\frak p\ar[r]\ar@{=}[d]& G(R/\frak p)\ar[d]\ar[r]&I\ar[d] \ar[r]& 0\\
0\ar[r]&R/\frak p\ar[r]&E(R/\frak p) \ar[r]& G\ar[r]& 0}$$
where $I$ is injective as $\Gid R/\frak p=1$ and we notice that $I_{\frak p}=0$ as $k(\frak p)\cong G(k(\frak p))=G(R/\frak p)_{\frak p}$. Moreover, the diagram gives rise to the exact sequence $0\To G(R/\frak p)\To E(R/\frak p)\oplus I\To G\To 0.$
 Localizing at $\frak p$ and  the above arguments give rise to the exact sequence of $R_{\frak p}$-modules $0\To k(\frak p)\To E(R/\frak p)\To k(\frak p)\To 0$
 and application of the functor $\Hom_{R_{\frak p}}(k(\frak p),-)$ gives a exact sequence of $k(\frak p)$-vector spaces $0\To \Hom_{R_{\frak p}}(k(\frak p),k(\frak p))\To k(\frak p)\To \Hom_{R_{\frak p}}(k(\frak p),k(\frak p))\To \Ext_{R_{\frak p}}^1(k(\frak p),k(\frak p))\To 0$ and the isomorphism  $\Ext_{R_{\frak p}}^i(k(\frak p),k(\frak p))\cong\Ext_{R_{\frak p}}^{i+1}(k(\frak p),k(\frak p))$ which yield $\Ext_{R_{\frak p}}^i(k(\frak p),G(R/\frak p)_{\frak p})\cong \Ext_{R_{\frak p}}^i(k(\frak p),k(\frak p))\cong k(\frak p)$ for each $i\geq 0$. 
  \end{Example}
\medskip

\begin{Theorem}\label{bass}
Let $\frak p\in\Spec R$ such that $\Gpd_RR/\frak p$ is finite. Then $\mu^i(\frak p,G(R/\frak p))$ is finite for each $i\geq 0$. 
\end{Theorem}
\begin{proof}
It follows from [CFH, Theorem 5.3] that $\Gpd_{R_{\frak p}}(R_{\frak p}/\frak pR_{\frak p})$ is finite.
By \cref{corlocc}, we have $G(R/\frak p)_{\frak p}=G(R_{\frak p}/\frak pR_{\frak p})$ and since $\mu^i(\frak p,G(R/\frak p))=\mu^i(\frak pR_{\frak p},G(R_{\frak p}/\frak pR_{\frak p}))$ for each $i\geq 0$, we may assume that $R$ is a local ring with the maximal ideal $\frak m$ and the residue field $k=R/\frak m$. By virtue of [H1, Corollary 2.13], there exists an exact sequence of $R$-modules $0\To L\To G\To k\To 0$ such that $G$ is a finitely generated Gorenstein projective $R$-module and $L$ has finite projective dimension. Application of the functor $(-\breve{)}=\Hom_R(-,E(k))$, we have an exact sequence of $R$-modules $0\To k\To \breve{G}\To \breve{L}\To 0.$ We notice that $\breve{M}$ is Gorenstein injective and $\breve{L}$ has finite injective dimension and so $k\To \breve{G}$ is a Gorenstein injective preenvelope of $k$. Therefore $G(k)$ is a direct summand of $\breve{G}$ and since $\breve{G}$ is artinian, we deduce that $G(k)$ is artinian; and hence $\mu^i(\frak m,G(k))$ is finite for each $i\geq 0$. 
\end{proof}

\medskip
\begin{Proposition}\label{two}
Let $\frak p\in\Spec R$ such that $R_{\frak p}$ is a Gorenstein local of dimension one. Then $G(k(\frak p))/k(\frak p)\cong E(R/\frak p)$, where $k(\frak p)=R_{\frak p}/\frak pR_{\frak p}$. 
\end{Proposition}
\begin{proof}
 There is an  exact sequence $0\To R/\frak p\To G(R/\frak p)\To E\To 0$ such that $E\in ^{\bot}\cG\cI$. By [CFH, Theorem 6.3], we have $\Gid_{R_{\frak p}}k(\frak p)=1$ and so by \cref{tor}, 
 we have $E_{\frak p}=E(R/\frak p)^{(Y)}$. Since $\depth R_{\frak p}=1$, there exist $x\in\frak pR_{\frak p}$ and an exact sequence of $R_{\frak p}$-modules $0\To R_{\frak p}\stackrel{x.}\To R_{\frak p}\To X\To 0$ where $X=R_{\frak p}/xR_{\frak p}$. Since $\depth X=0$, it contains a copy of $k(\frak p)$. Applying $\Hom_{R_{\frak p}}(-,k(\frak p))$, we deduce that $\Hom_{R_{\frak p}}(X,k(\frak p))\cong k(\frak p)\cong\Ext^1_{R_{\frak p}}(X,k(\frak p))$. On the other hand, there exists the following commutative diagram of $R_{\frak p}$-modules with the exact row
 
 $$\xymatrix{\Hom_{R_{\frak p}}(X,k(\frak p))\ar[r]&\Hom_{R_{\frak p}}(X,G(R/\frak p)_{\frak p})\ar[r]&\Hom_{R_{\frak p}}(X,E_{\frak p})\ar[r]\ar@{->>}[rd]& \Ext_{R_{\frak p}}^1(X,k(\frak p))\ar@{-->>}[d]^{\theta}\ar[r]& 0\\
&&&\Hom_{R_{\frak p}}(k(\frak p),E_{\frak p})}$$
where the zero map 
$\Hom_{R_{\frak p}}(X,G(R/\frak p)_{\frak p})\To\Hom_{R_{\frak p}}(X,E_{\frak p})\To \Hom_{R_{\frak p}}(k(\frak p),E_{\frak p})$ ensures the existence of the epimorphism $\theta$ which implies that $\dim_{k(\frak p)}\Hom_{R_{\frak p}}(k(\frak p),E_{\frak p})=1$.
\end{proof}

\section{Endomorphism rings of indecomposable Gorenstein injective modules}
 Assume that $\frak p\in\Spec R$ and set $S(\frak p)=\{f\in{\rm End}_R(G(R/\frak p))|\hspace{0.1cm} f(R/\frak p)\subset R/\frak p\}$. Clearly, $S(\frak p)$ is a ring (may be non-commutative). We have the following theorem.
\medskip

\begin{Theorem}
Let $\frak p\in\Spec R$. Then $S(\frak p)$ is a local ring containing $R/\frak p$ with the maximal ideal $\frak n=\{f\in S(\frak p)|\hspace{0.1cm} f(R/\frak p)=0\}$.
\end{Theorem}
\begin{proof}
It is clear that $\frak n$ is an ideal of $S(\frak p)$. We show that any non-unit element of $S(\frak p)$ is in $\frak n$. If $f\in S(\frak p)$ is non-unit and $f\notin\frak n$, then $f:R/\frak p\To R/\frak p$ is non-zero and clearly, it is injective.  Then we have the following commutative diagram 
$$\xymatrix{R/\frak p\ar[r]^{\cong}\ar[d]& \Im f \ar[d]\\
G(R/\frak p)\ar[r]^f&G(R/\frak p)}.$$
 Since $G(R/\frak p)$ is Gorenstein injective envelope of $R/\frak p\cong\Im f$, the map $f$ is an isomorphism which is a contradiction.  Since Rad$(S(\frak p))\subset\{$non-units$\}\subset\frak n$, it suffices to show that $\frak n\subset$Rad$(S(\frak p))$. Let $f\in\frak n$ and so $f(R/\frak p)=0$. Since $R$ admits a dualizing complex, $\bigoplus_{\bbN}G(R/\frak p)$ is Gorenstein injective  and by [EJ2, Corollary 6.4.2], it is the Gorenstein injective envelope of $\bigoplus_{\bbN} R/\frak p$. Taking the homomorphism $\theta:\bigoplus_{\bbN} G(R/\frak p)\To  \bigoplus_{\bbN}G(R/\frak p)$ by $\theta(x_0,x_1,\dots)=(x_0,x_1-f(x_0),x_2-f(x_1),\dots)$, we have $\theta|_{\bigoplus_{\bbN} R/\frak p}=1_{\bigoplus_{\bbN} R/\frak p}$ so that $\theta$ is an automorphism. Now for each $x\in G(R/\frak p)$, there exists $(x,x_1,\dots)\in\bigoplus_{\bbN}G(R/\frak p)$ and so $x_n=f^n(x)$ and since $x_n=0$ for all $n\gg 0$, there exists $n\in\bbN$ such that $f^n(x)=0$. Therefore $\Sigma_{n\geq 0} (-f)^n\in\End_R(G(R/\frak p))$ and $(1+f)(\Sigma_{n\geq 0} (-f)^n)=1_G(R/\frak p)$. 
 To prove the first assertion we define $\alpha:R/\frak p\To S$ by $\theta(x+p)=x.$ where $G(R/\frak p)\stackrel{x.}\To G(R/\frak p)$ is multiplication map. It is clear that $\alpha$ is a ring homomorphism and since $\frak p\in\Ass_RG(R/\frak p)$, the map $\alpha$ is injective. 
\end{proof}

 \medskip
 \begin{Lemma}\label{free}
 Let $\frak p\in\Spec R$ and $n\geq 1$. Then $\End_R(E(R/\frak p)^n)\cong \widehat{R_{\frak p}}^{n^2}$. In particular,
 $\End_R(E(R/\frak p)^n)$ is a notherian $\widehat{R_{\frak p}}$-module. 
 \end{Lemma}
 \begin{proof}
 We have the following isomorphisms 
 $$\End_{R}(E(R/{\frak p})^n)\cong\Hom_{R_{\frak p}}(E(R/\frak p),E(R/\frak p)^n)^{n}$$$$\cong (\underset{\leftarrow}{\rm lim}\Hom_{R_{\frak p}}(\Hom_{R_{\frak p}}(R_{\frak p}/\frak p^nR_{\frak p},E(R/\frak p)), E(R/\frak p))^{n}))^{n}$$$$\cong (\underset{\leftarrow}{\rm lim}(\Hom_{R_{\frak p}}(\Hom_{R_{\frak p}}(R_{\frak p}/\frak p^tR_{\frak p},E(R/\frak p)), E(R/\frak p))^{n}))^{n}\cong(\underset{\leftarrow}{\rm lim}(R_{\frak p}/\frak p^tR_{\frak p})^{n})^{n}\cong\widehat{R_{\frak p}}^{n^2}$$
 where the third and fourth isomorphism hold because the $R_{\frak p}$-modules $\Hom_{R_{\frak p}}(R_{\frak p}/\frak p^tR_{\frak p},E(R/\frak p))$  and $R_{\frak p}/\frak p^tR_{\frak p}$ have finite length and the last isomorphism holds as $R_{\frak p}^{n}$ is a finitely generated $R_{\frak p}$-module.  This implies that $\End_{R}(E(R/{\frak p})^n)$ is a free $\widehat{R_{\frak p}}$-module.
 \end{proof}
 
 \medskip

 In the rest of this section, for any $\frak p\in\Spec R$, we set $S(\frak p)_{\frak p}=\{f_{\frak p}\in\End_{R_{\frak p}}(G(k(\frak p))|\hspace{0.1cm}f\in S(\frak p)\}$. We now the following lemma.
  
 \medskip
 \begin{Lemma}
  Let $\frak p\in\Spec R$. Then $S(\frak p)_{\frak p}$ is local. 
  \end{Lemma}
  \begin{proof}
  We claim that $\frak n_{\frak p}$ is the maximal ideal of $S(\frak p)_{\frak p}$. Given an arbitrary proper ideal $I$ of $S(\frak p)_{\frak p}$, consider $\frak a=\{f\in S(\frak p) |\hspace{0.1cm} f_{\frak p}\in I\}$. It is clear that $\frak a$ is an ideal of $S(\frak a)$ and so $\frak a\subset\frak n$. If $1_{G(R/\frak p)_{\frak p}}\in \frak n_{\frak p}$, then $k(\frak p)=0$ which is a contradiction; consequently $\frak n_{\frak p}$ is a proper ideal of $S(\frak p)_{\frak p}$.
   \end{proof}

\medskip
\begin{Theorem}\label{ring}
Let $\frak p\in\Spec R$ such that $\mu(\frak p,G(R/\frak p))$ is finite. Then $S(\frak p)_{\frak p}$ is noetherian. 
\end{Theorem}
\begin{proof}
 Assume that $\mu(\frak p,G(R/\frak p))=n$. The fact that $G(R/\frak p)$ is a $\frak p$-torsion module and the assumption implies that $G(R/\frak p)_{\frak p}=G(k(\frak p))$ is a $\widehat{R_{\frak p}}$-submodule of  $E(R/\frak p)^n$. On the other hand, since $\Hom_R(E(R/\frak p),E(R/\frak p))\cong \widehat{R_{\frak p}}$, we have $\End_R(E(R/\frak p)^n)\cong\bbM_n(\widehat{R_{\frak p}})$; where $\bbM_n(\widehat{R_{\frak p}})$ is the ring of $n\times n$-matrices over $\widehat{R_{\frak p}}$.
 We first prove that $\bbM_n(\widehat{R_{\frak p}})$  is a noetherian ring containing $\widehat{R_{\frak p}}$ and so is $\End_R(E(R/\frak p)^n)$. If $I_1\subset I_2\subset\dots $ is a ascending chain of ideals of $\bbM_n(\widehat{R_{\frak p}})$, by [L, Corollary 17.8], we have $I_i=\bbM(\frak a_i)$ for some ideal $\frak a_i$ of $\widehat{R_{\frak p}}$. Assume that $E_{11}=(e_{ij})$ be a $n\times n$-matrices with $e_{ij}=0$ for all $i,j\neq 1$, but $e_{11}=1$. Then for every $x\in\frak a_i$, we have $xE_{11}\in I_i$ and so $xE_{11}\in I_{i+1}=\bbM(\frak a_{i+1})$ so that $x\in\frak a_{i+1}$. Therefore $\frak a_1\subset \frak a_2\subset\dots $ is an ascending chain of ideals of noetherian ring $\widehat{R_{\frak p}}$. Then there exists an integer $t$ such that $\frak a_t=\frak a_i$ for all $i>t$; consequently $I_t=I_i$ for all $i>t$. The rings homomorphism $\theta:\widehat{R_{\frak p}}\To \bbM(\widehat{R_{\frak p}})$; given by $x\mapsto xI_n$ is injective where $I_n=(\delta_{ij})=1_{\bbM(\widehat{R_{\frak p}}))}$ so that $\widehat{R_{\frak p}}$ is a subring of $\bbM(\widehat{R_{\frak p}}))$. 
 We now consider $A=\{f\in\End_R(E(R/\frak p)^n)|\hspace{0.1cm}f|_{G(k(\frak p))}\in S(\frak p)_{\frak p}\}$; where $f|_{G(k(\frak p))}$ means the restriction of $f$ to $G(k(\frak p))$. It is clear that $A$ is a subring of $\End_R(E(R/\frak p)^n)$. On the other hand, since by \cref{tor}, the module $G(k(\frak p))$ is $\frak p$-torsion, for every $g\in G(k(\frak p))$, there exists a positive integer $n$ such that $\frak p^ng=0$. Then for any $(r_i)\in\widehat{R_{\frak p}}$, the multiplication $(r_i)g:=r_ng$ arranges  $G(k(\frak p))$ as an $\widehat{R_{\frak p}}$-module. Thus the rings homomorphism $\alpha:\widehat{R_{\frak p}}\To A$, by $r\mapsto r.$ implies that $\widehat{R_{\frak p}}$ is a subring of $A$. We observe that $A$ is an $\widehat{R_{\frak p}}$-submodule of the $\widehat{R_{\frak p}}$-module $\End_R(E(R/\frak p)^n)$ and also any right and left ideal of $A$ is an $\widehat{R_{\frak p}}$-submodule of the $\widehat{R_{\frak p}}$-module $\End_R(E(R/\frak p)^n)$ and so \cref{free} implies that $A$ is a noetherian ring. Now, the rings homomorphism $\theta:A\To S(\frak p)_{\frak p}$, given by $f\mapsto f|_{G(k(\frak p))}$ is surjetive because any $g\in S(\frak p)_{\frak p}$ can be extended to a homomorphism $f\in A$. Consequently, $S(\frak p)_{\frak p}$ is noetherian. 
 \end{proof}

\medskip
\begin{Corollary}
Let $R$ be a Gorenstein ring of Krull dimension $d$. Then $S(\frak p)_{\frak p}$ is noetherian all $\frak p\in\Spec R$. 
\end{Corollary}
\begin{proof}
Straightforward using \cref{goro} and \cref{ring}. 
\end{proof}

\medskip

\begin{Corollary}
Let $\frak p\in\Spec R$ such that $\Gid_RR/\frak p$ is finite. Then $S(\frak q)_{\frak q}$ is noetherian for all $\frak q\in\Spec R$ with $\frak q\subset \frak p$.
\end{Corollary}
\begin{proof}
Straightforward using \cref{subb} and \cref{ring}. 
\end{proof}

\medskip

\begin{Proposition}\label{mat}
Let $\frak p\in\Spec R$. Then $\widehat{R_{\frak p}}/\Ann_{\widehat{R_{\frak p}}}G(R/\frak p)_{\frak p}$  is a (commutativ) subring of $S(\frak p)_{\frak p}$. 
\end{Proposition}
\begin{proof}
There is a  rings homomorphism $\alpha: \widehat{R_{\frak p}}\To S(\frak p)_{\frak p}$, given by $r\mapsto r.$. Clearly $\Ker \alpha=\Ann_{\widehat{R_{\frak p}}}G(R/\frak p)_{\frak p}$.
\end{proof}

\medskip

\begin{Corollary}
 Let $\frak p\in\Spec R$ such that $R_{\frak p}$ is a Gorenstein local ring of dimension one. Then $S(\frak p)_{\frak p}\cong\widehat{R_{\frak p}}$.
\end{Corollary}
\begin{proof}
By [CFH, Theorem 6.3], we have $\Gid_{R_{\frak p}}(k(\frak p)=1$ and so in the exact sequence $0\To R/\frak p\To G(R/\frak p)\To E\To 0$, the module $X$ is injective and it follows from \cref{two} that $E_{\frak p}=E(R/\frak p)$ and so there exists a  rings homomorphism $\theta:S(\frak p)_{\frak p}\To \widehat{R_{\frak p}}$; given by $f_{\frak p}\To r$ satisfying the following commutative diagram of $\widehat{R_{\frak p}}$-modules 

$$\xymatrix{0\ar[r]&k(\frak p)\ar[r]\ar[d]& G(k(\frak p))\ar[d]^{f_{\frak p}}\ar[r]&E(R/{\frak p})\ar[d]^{r.} \ar[r]& 0\\
0\ar[r]&k(\frak p)\ar[r]& G(k(\frak p))\ar[r]& E(R/{\frak p})\ar[r]& 0.}$$
 
  For every $r\in\widehat{R_{\frak p}}$, we have $\theta(r.)=r$ and so $\theta$ is surjective. For injectivity $\theta$, if $\theta(f_{\frak p})=0$ and $f_{\frak p}$ is non-zero, then $f_{\frak p}(G(k(\frak p)))\subset k(\frak p)$ and so $f_{\frak p}:G(k(\frak p))\To k(\frak p)$ is surjective. Since $\depth R_{\frak p}\geq 1$, there exists a non-zero-divisor $x\in\frak p$ and so the exact sequence $0\To R_{\frak p}\stackrel{x.}\To R_{\frak p}\To R_{\frak p}/xR_{\frak p}\To 0$ implies that $xG(k(\frak p))=G(k(\frak p))$ so that $xk(\frak p)=f(xG(k(\frak p)))=k(\frak p)$. But this implies that $k(\frak p)=0$; which is a contradiction.  
\end{proof}


\end{document}